# Two methods for solving optimization problems arising in electronic measurements and electrical engineering

Ya.D. Sergeyev*    P. Daponte†    D. Grimaldi‡    A. Molinaro‡

March 19, 2011

**Abstract**

In this paper we introduce a common problem in electronic measurements and electrical engineering: finding the first root from the left of an equation in the presence of some initial conditions. We present examples of electrotechnical devices (analog signal filtering), where it is necessary to solve it. Two new methods for solving this problem, based on global optimization ideas, are introduced. The first uses the exact a priori given global Lipschitz constant for the first derivative. The second method adaptively estimates local Lipschitz constants during the search. Both algorithms either find the first root from the left or determine the global minimizers (in the case when the objective function has no roots). Sufficient conditions for convergence of the new methods to the desired solution are established in both cases. The results of numerical experiments for real problems and a set of test functions are also presented.

## 1   Introduction

Let us consider a device whose behaviour depends on a characteristic $f(x)$, $x \in [a, b]$; $f(x)$ may be, for instance, an electrical signal and $[a, b]$ a time interval. The device functions correctly while $f(x) > 0$. Of course, at the initial moment $x = a$ we have $f(a) > 0$. It is supposed that the function $f(x)$ is multiextremal and the Lipschitz condition is satisfied for its first derivative, i.e.

*ISI-CNR c/o DEIS, Universita' della Calabria, 87036 Rende (CS), Italy and Nizhni Novgorod State University, pr. Gagarina 23, Nizhni Novgorod, Russia, e-mail: yaro@si.deis.unical.it

†Dipartimento di Ingegneria dell'Informazione ed Ingegneria Elettrica, via Ponte Don Melillo 1, 84084 Fisciano (SA) - Italy, e-mail: daponte@nadis.dis.unina.it

‡Dipartimento di Elettronica, Informatica e Sistemistica, Universita' della Calabria, 87036 Rende (CS) - Italy, e-mail: grimaldi@ccusc1.unical.it

$$|f'(x) - f'(y)| \leq L|x-y|, \quad x,y \in [a,b] \qquad (1)$$

where the constant $L$, $0 < L < \infty$, is the Lipschitz constant. Generally $L$ is unknown. The problem we deal with in this paper is determining a time interval $[a, x^*]$ where the device works correctly. This problem is equivalent to the problem of finding the root of $f$ which is first from the left, that is, finding

$$x^* = min\{x : f(x) = 0\}, \qquad (2)$$

subject to conditions:

$$x \in [a,b], f(a) > 0. \qquad (3)$$

This problem very often arises in electronic measurements ([1], [6], [23]) and electrical engineering ([3], [5], [11]–[14], [18]) and in Section 2 we present two concrete applications as examples. Usually it is difficult to solve the problem (2), (3) in an analytical way and numerical methods are used to find a $\sigma$-approximation $x_\sigma$ of the point $x^*$ such that

$$0 \leq f(x_\sigma), |x_\sigma - x^*| < \sigma. \qquad (4)$$

Two approaches are currently used by engineers for solving the problem (1)–(4). The first one uses standard local techniques for finding equation roots in order to achieve a rapid convergence to the point $x^*$. The drawback of these methods is that convergence is not assured since $f(x)$ is a multiextremal function on $[a,b]$ and the methods may diverge or converge to a local minimum greater than zero (see [17]). Moreover, if the objective function $f(x)$ has more than one root (and this is usually the case, see Fig.1), different choices of the initial conditions can produce different solutions of the equation $f(x) = 0$.

The second approach is based on the use of any simple grid technique which produces a dense mesh starting from the left margin of the interval and going on by $\sigma$ till the value $f(x)$ becomes less than zero (see [3]). This approach is very reliable but the number of evaluations of $f(x)$ is too high.

In this paper we propose two new numerical algorithms for solving the problem (1)–(4) in order to find a point $x_\sigma$ from (4). The methods are based on geometric ideas of the global optimization technique [21]. The new algorithms either find the point $x_\sigma$ from (4) or determine a $\sigma$ approximation $x'_\sigma$ of the global minimizer $x'$ of $f(x)$ and the corresponding value $f(x'_\sigma)$ in the case

$$f(x) > 0, x \in [a,b]. \qquad (5)$$

The first algorithm uses the given constant $L$ from (2). Since in practice it is difficult to know this value a priori, the problem of estimating it in the course of the search arises. The second method presented here estimates the local Lipschitz constants for subintervals of the search region in the course of the search. Using local estimates instead of the global accelerates the search significantly.

The new techniques are described in Section 3. Convergence conditions are established in Section 4. Section 5 contains results of numerical experiments executed both with objective functions derived from the applications presented in Section 2 and with a series of test functions. Finally, Section 6 concludes the paper.

## 2 Filters as examples of electronic devices where the problem arises

The problem (1)–(4) arises very often in applications (see [1], [3], [5], [6], [11]–[14], [18], [23]). In fact, the objective function $f(x)$ can be considered as, for example, a reliability characteristic of a device or a mathematical model. While $f(x) > 0$ the model is reliable and for $x \geq x^*$ it is not.

Here we consider two concrete examples leading to illustrate the problem (1)–(4). Both of them deal with electrical filters. *Filters* are basic electronic components utilized in many fields such as power conversion circuits, electronic measurement instruments and communications systems. In particular, *electrical filters* can be found in the telephone, television, radio, radar and sonar. A filter is a device that modifies in a predetermined way the input signal that passes through it. Electrical filters may be classified as: *analog* filters, used to process analog or continuous-time signals; *digital* filters, used to process digital signals (discrete-time signals).

Let us consider a signal $s(x)$, where $x$ is time. If a signal $s(x)$, composed of a sum of signals $s_1(x), s_2(x), ..., s_n(x)$ so that

$$s(x) = s_1(x) + s_2(x) + .. + s_n(x),$$

is the input of an analog filter, the output signal is obtained from the input one by suppressing certain components $s_k(x), k \in \{1, .., n\}$. Let us define for the signal $s(x)$ its *frequency* $\theta$ as the number of times that the signal repeats itself in unit time and the *pulse* $\omega = 2\pi\theta$. Below we refer to $\theta$ or $\omega$ simply as frequency. If a signal $s(x)$ has a certain frequency $\theta$, it may be represented by a rotant vector having angular speed equal to $\omega$ and the amplitude equal to the maximum amplitude of $s(x)$. As all the vectors with the same angular speed can be represented in a complex plane, (see [5]) since we can represent the signal $s(x)$ in the frequency domain instead of the time domain, so that we can write $s(x)$ as $s(\omega)$. If $Y(\omega)$ is the filter output and $X(\omega)$ is its input into the frequency domain, the function:

$$|H(\omega)| = \frac{|Y(\omega)|}{|X(\omega)|}$$

is called the *transfer function* in the frequency domain (see [11], [12]). From the value of $|H(\omega)|$ we can evaluate the answer of the filter, that is, how far the output signal is from the input signal. The *cutoff frequency* $\omega_c$ is defined as the frequency where the transfer function is equal to $\sqrt{0.5}$ of its maximum amplitude. This implies that $\omega_c$ can be calculated from the following equation:

$$|H(\omega)| = \sqrt{0.5} H_{max},$$

or

$$|H(\omega)|^2 = 0.5 H_{max}^2$$

where

$$H_{max} = max\{H(\omega) : \omega \in [0, \infty)\}.$$

The *passband* is the width of the interval $[\omega_{c1}, \omega_{c2}]$ in which

$$|H(\omega)|^2 \geq 0.5 H_{max}^2,$$

where $\omega_{c1}$ is called the *lower cutoff frequency* and $\omega_{c2}$ is called the *higher cutoff frequency*. If $\omega_{c1} = 0$ the filter is a *low-pass filter*; if $\omega_{c2} \to \infty$ the filter is a *high-pass filter*; finally, if $\omega_{c1}$ and $\omega_{c2}$ are finite values, the filter is a *passband filter*. In general, an electrical filter is designed to separate one component of the input signal from the others. If we want to separate a particular frequency $\omega_p$ and reject all other frequencies, for technical reasons we cannot build a filter that will pass only $\omega_p$, but a set of frequencies $\omega_p \in [\omega_{c1}, \omega_{c2}]$. As an example, let us consider a radio or a television receiver. The transmission station is assigned an interval of frequencies called the *band of frequencies* or *channel frequencies*, in which it must transmit its signal. Ideally, the receiver should accept and process any signal in the assigned channel and completely exclude signals at all other frequencies. Thus the simplest specifications on the magnitude of the transfer function of the receiver are:

$$|H(\omega)| = \begin{cases} H_{max} & \text{for } \omega_{c1} \leq \omega \leq \omega_{c2} \\ 0 & \text{otherwise} \end{cases} \quad (6)$$

where $\omega_{c1} \leq \omega \leq \omega_{c2}$ is the channel of the signal to be received. However, no circuits can produce such a transfer function exactly. In practice, filters are not required to meet the extremely stringent requirements such as those of (6) and some filters with a transfer function approximating (6) have been found to be consistently satisfactory. A filter that has a uniform approximation property in the passband is the Chebyshev filter (see [11], [12]), our first example that will be described in the following.

An electrical circuit which realizes a Chebyshev filter is shown in Fig.2, where voltage $V_{in}(\omega)$ is the input and voltage $V_{out}(\omega)$ is output. The transfer function $F(\omega)$, obtained applying Kirchhoff laws to the circuit of Fig.2, has the following expression:

$$F(\omega) = |\frac{V_{out}(\omega)}{V_{in}(\omega)}| = \frac{1}{\sqrt{1 + R^2 C^2 \omega^2}} \cdot \frac{1}{\sqrt{(2 - \omega^2 LC)^2 + \omega^2 L^2/R^2}}, \quad (7)$$

where the values $R$, $C$ and $L$ are introduced in Fig.2 and $|\cdot|$ is the length of a complex vector. This function suppresses the frequency components beyond the cutoff frequency $\omega_c$ and is characterized by ripples in the passband. The number of maxima and minima in the ripple passband defines the filter order. In our case, the filter order is $n = 3$. The cutoff frequency can be found as the first root from the left for the function

$$f(\omega) = F(\omega)^2 - 0.5F_{max}^2,$$

where $F_{max}$ is the maximum of the function

$$F(\omega), \omega \in (0, \infty).$$

Let us consider the second example. In Fig.3 the electrical circuit of a bandpass filter [11] is shown. The transfer function of this filter is given by:

$$F(\omega) = |\frac{V_{out}(\omega)}{I_{in}(\omega)}| = \frac{\omega L_1 R_1}{\sqrt{(Z_1^2 + Z_2^2)^2 \cdot Z_3}} \qquad (8)$$

where:

$$Z_1 = -\omega^3 R_1 L_1 L_2 + \omega R_1 L_2 + \omega R_1 L_1 C_1/C_2 - R_1/(\omega C_2) + 2\omega L_1 R_1 + \omega L_1 R_2$$

$$Z_2 = \omega^2 L_1 L_2 + \omega^2 R_1 R_2 L_1 C_1 - R_1 R_2 - L_1/C_2$$

$$Z_3 = (\omega L_1)^2 + (\omega^2 R_1 L_1 C_1 - R_1)^2$$

This result is obtained applying Kirchhoff laws to the circuit in Fig.3. The transfer function tends to zero for $\omega \to 0$ and $\omega \to \infty$. The cutoff frequency can be found as the first zero crossing for the function

$$f(\omega) = -(F(\omega)^2 - 0.5F_{max}^2).$$

The values of the circuit parameters may be changed, so varying the flatness of the transfer function.

## 3 Theoretical background and the algorithms

New algorithms presented here for solving the problem (1)–(4) are based on the following idea. Let us suppose that the objective function $f(x)$ and its first derivative $f'(x)$ have been already calculated at $n$ trial points $x^i, 1 \leq i \leq n$. We can reorder these points by subscripts in such a way that

$$a = x_1 < x_2 < ... < x_i < ... < x_n = b.$$

Thus, dealing with the trial points, we use two numerations. The record $x^i$ means that this point has been produced during the $i$th iteration. The record $x_i$, $i = i(n)$, means that this point has the $i$th position in the string above during the $n$th iteration. We designate the results of trials as $z_i = f(x_i)$, $z'_i = f'(x_i)$, $1 \leq i \leq n$.

For every interval $[x_{i-1}, x_i]$, $1 < i \leq n$, we construct an auxiliary function $\phi_i(x)$ in such a way that $\phi_i(x) \leq f(x)$, $x \in [x_{i-1}, x_i]$. Knowing the structure of the auxiliary function, we can find the first function $\phi_i(x)$ such that for some $x \in [x_{i-1}, x_i]$ it follows

$$\phi_i(x) = 0$$

and determine the first root from the left of this equation. Adaptively improving the set of functions $\phi_i(x)$, $1 < i \leq n$, by adding new points $x^{n+1}, x^{n+2}, ...$ we will improve our approximation of $f(x)$ and of the point $x^*$. This geometric approach is widely used in global optimization (see [7], [10]) applying functions $\phi_i(x)$ with different structures (see [8], [15], [16] for methods using only the values of objective functions and [2], [4], [9], [21], [25], [26], [27], [28] for algorithms where the first derivatives are also taken into consideration).

In the two algorithms presented here we use the following three ideas to provide a fast localization of the points $x_\sigma$ from (4):

- smooth auxiliary functions from [21] where they demonstrated good performance in global optimization;

- constructing auxiliary functions only for intervals $[x_{i-1}, x_i]$, $1 < i \leq k$, where

$$k = min\{\{n\} \cup \{i : f(x_i) < 0, 1 < i \leq n\}\}. \tag{9}$$

- adaptive estimating of the *local* Lipschitz constant $L_i$ for every interval $[x_{i-1}, x_i]$ (in the second algorithm).

Let us discuss these ideas one after another. First, there exist three types of support function elaborated to solve global optimization problems: piece-wise linear (see [8], [15], [16], [19]), nonsmooth quadratic (see [4], [9], [25], [27]), and smooth quadratic (see [21], [26], [28]). We use the last construction because (see [21], [26], [28]) it best approximates the objective function.

Suppose that we have an estimate $m_i$ of the constant $L_i$ such that:

$$m_i \geq L_i. \tag{10}$$

In this case it is possible to construct a support function $\phi_i(x)$ for $f(x)$ over $[x_{i-1}, x_i]$ (see [21]) as follows:

$$\phi_i(x) = \begin{cases} z_{i-1} + z'_{i-1}(x - x_{i-1}) - 0.5m_i(x - x_{i-1})^2, & x \in [x_{i-1}, y'_i] \\ 0.5m_i x^2 + b_i x + c_i, & x \in (y'_i, y_i] \\ z_i - z'_i(x_i - x) - 0.5m_i(x_i - x)^2, & x \in (y_i, x_i] \end{cases} \tag{11}$$

where

$$y_i = \frac{x_i - x_{i-1}}{4} + \frac{z'_i - z'_{i-1}}{4m_i} + \frac{z_{i-1} - z_i + z'_i x_i - z'_{i-1} x_{i-1} + 0.5 m_i (x_i^2 - x_{i-1}^2)}{m_i (x_i - x_{i-1}) + z'_i - z'_{i-1}}, \quad (12)$$

$$y'_i = -\frac{x_i - x_{i-1}}{4} - \frac{z'_i - z'_{i-1}}{4m_i} + \frac{z_{i-1} - z_i + z'_i x_i - z'_{i-1} x_{i-1} + 0.5 m_i (x_i^2 - x_{i-1}^2)}{m_i (x_i - x_{i-1}) + z'_i - z'_{i-1}}, \quad (13)$$

$$b_i = z'_i - 2 m_i y_i + m_i x_i, \quad (14)$$

$$c_i = z_i - z'_i x_i - 0.5 m_i x_i^2 + m_i y_i^2. \quad (15)$$

An illustration of the functions $f(x)$, $\phi_i(x)$ is presented in Fig.1. The function $\phi_i(x)$ has been constructed using the Taylor formula for $f(x)$ about the point $x_{i-1}$ (see the first line in (11)) and the point $x_i$ (see the third line in (11)). The second line of (11) has been obtained as the quadratic which agrees with the $f(x)$ curvature at the interval extrems. Note that the first derivative $\phi'_i(x)$, for all $x \in (x_{i-1}, x_i)$ exists.

It will be useful for us to find the point

$$h_i = argmin\{\phi_i(x) : x \in [x_{i-1}, x_i]\} \quad (16)$$

and the corresponding value

$$R_i = \phi_i(h_i) = min\{\phi_i(x) : x \in [x_{i-1}, x_i]\}. \quad (17)$$

We will call the value $R_i$ the *characteristic* of the interval $[x_{i-1}, x_i]$. Let us consider two cases. If $\phi'_i(y'_i) < 0$ and $\phi'_i(y_i) > 0$, then

$$h_i = argmin\{f(x_{i-1}), \phi_i(\widehat{x}_i), f(x_i)\} \quad (18)$$

where:

$$\widehat{x}_i = 2 y_i - z'_i m_i^{-1} - x_{i-1}. \quad (19)$$

The point $\widehat{x}_i$ is determined from the equation $\phi'_i(x) = 0$, $x \in [y'_i, y_i]$. It follows from (11) that

$$\phi_i(\widehat{x}_i) = c_i - 0.5 m_i \widehat{x}_i^2. \quad (20)$$

In the second case there is no point $\widehat{x}_i \in [y'_i, y_i]$ such that $\phi'_i(\widehat{x}_i) = 0$ and

$$h_i = argmin\{f(x_{i-1}), f(x_i)\}. \quad (21)$$

The algorithms work by constructing the function $\phi_i(x)$ from (11) from left to right taking the intervals one after another and calculating their characteristics. If in a step $R_j < 0$ has been found, this means that there exists a point $\tilde{x} \in [x_{j-1}, x_j]$ such that $\phi_j(\tilde{x}) = 0$.

In this case we determine the new trial point $x^{n+1} = \tilde{x}$ (all possible locations of $x^{n+1}$ are shown in Figs. 4-6) and evaluate $f(x^{n+1})$ and $f'(x^{n+1})$. If $f(x^{n+1}) < 0$ then there is no need to consider the interval $(x^{n+1}, b]$ because the solution $x_\sigma$ is in $(a, x^{n+1}]$ (here $a$, $b$ are from (3)). Then we increment $n$ and restart the procedure. If (5) takes place, then the algorithm will function to find an approximation $x'_\sigma$ of the point $x'$.

The last point we discuss before describing the algorithms are the obtaining of values $m_i$ such that (10) holds and the influence on the speed (and correctness) of the algorithm. The first way is to take $m_i = L$, where $L$ is from (1). But in real problems often it is difficult to know the exact value of $L$. In this case a fixed estimate of $L$ is taken and used in the course of the search. This strategy is applied in global optimization methods (for different auxiliary functions $\phi_i(x)$) in [2], [4], [15], [28]. We use it in our first algorithm (A1).

The second approach to determining the values $m_i$ is to estimate $L$ in the course of the search using information obtained from evaluating $f(x)$, $f'(x)$ at trial points. The way of estimating $L$ is under intensive investigation (see [24]) and many global optimization algorithms do it in different manners (see [8], [9], [16], [22]).

The main drawback of both approaches is the following. The global Lipschitz constant $L$ gives very poor information about the behaviour of $f(x)$ in every small interval $[x_{i-1}, x_i]$. That is why in our second algorithm (A2) we estimate local Lipschitz constants $L_i$ for every interval $[x_{i-1}, x_i]$, $1 < i \leq k$. This strategy has been succesfully applied in global optimization techniques [19]–[21], [25]–[27] for different classes of problems.

Now we are ready to describe the methods. We present only A2 since the A1 scheme can be easily obtained from it by eliminating **Step 2** and using $m_i = K$, $L \leq K < \infty$, in all subsequent steps.

Let us suppose that $n$ trials, with $n \geq 2$, of the algorithm have already been carried out at points $x^1, ..., x^n$. The $(n+1)$th trial point $x^{n+1}$ is chosen according to the following procedure.

**Step 1.** *(Ordering trial points)*

Among the trial points $x^1, ..., x^n$ of the previous $n$ iterations, form the subset $X^{k(n)}$ such that

$$X^{k(n)} = \{x_1, x_2, ..., x_k\},$$

where $k$ is defined in (9). We let $b^n = x_k$, the right margin of the search interval during the $n$th iteration (see Fig.1).

**Step 2.** *(Computing $m_i$)*

Calculate estimates $m_i$ for the local Lipschitz constants $L_i$ for the intervals

$$[x_{i-1}, x_i], 1 < i \leq k,$$

as follows:

$$m_i = r \cdot max\{\lambda_i, \gamma_i, \xi\}, \qquad (22)$$

where $\xi > 0$ and $r > 1$ are parameters of the method.

The values $\lambda_i$ and $\gamma_i$ spy on changes of local and global information respectively obtained in the course of the search. The value $\lambda_i$ is calculated as

$$\lambda_i = max\{v_j : 1 < j \le k, i-1 \le j \le i+1\}, \qquad (23)$$

where

$$v_j = \frac{|2(z_{j-1} - z_j) + (z'_j + z'_{j-1})(x_j - x_{j-1})| + d_j}{(x_j - x_{j-1})^2}, \qquad (24)$$

and

$$d_j = \sqrt{[2(z_{j-1} - z_j) + (z'_j - z'_{j-1})(x_j - x_{j-1})]^2 + (z'_j - z'_{j-1})^2(x_j - x_{j-1})^2}.$$

The second component $\gamma_i$ from (22) is calculated as

$$\gamma_i = m(x_i - x_{i-1})/X^{max} \qquad (25)$$

where $m$ estimates the global Lipschitz constant $L$ from (1)

$$m = max\{v_i : 1 < i \le k\} \qquad (26)$$

and

$$X^{max} = max\{x_i - x_{i-1} : 2 \le i \le k\}.$$

**Step 3.** *(Calculating characteristics $R_i$ of the intervals)*

Initialise the index sets $I = \emptyset$; $Y = \emptyset$; $Y' = \emptyset$. Set the index of the current interval $i = 2$.

**Step 3.0.** *(Organizing the main cycle)*

If $i > k$ then go to **Step 4**, otherwise compute the values $y_i$, $y'_i$ according to (12) and (13). If $\phi'_i(y'_i) \cdot \phi'_i(y_i) < 0$ then go to **Step 3.2**, otherwise go to **Step 3.1**.

**Step 3.1.** *(Computing $R_i$ if $\phi'_i(y'_i) \cdot \phi'_i(y_i) \ge 0$)*

Calculate $R_i = \phi_i(h_i)$, where $h_i$ is from (21). If $h_i = x_i$ then include $i$ in $Y$ else include $i$ in $Y'$. Go to **Step 3.3**.

**Step 3.2.** *(Computing $R_i$ if $\phi'_i(y'_i) \cdot \phi'_i(y_i) < 0$)*

Calculate $R_i = \phi_i(h_i)$, where $h_i$ is from (18). Include $i$ in $I$. Go to **Step 3.3**.

**Step 3.3.** *(Verifying the sign of $R_i$)*

If $R_i \le 0$ then go to **Step 5** otherwise set $i = i + 1$ and go to **Step 3.0**.

**Step 4.** *(Computing the new trial point if $R_j > 0, 1 < j \le k$)*

Find an interval $i$ with the minimal characteristic, i.e.

$$i = argmin\{R_j : 1 < j \le k\} \qquad (27)$$

and define the new trial at the point $x^{n+1}$ as follows

$$x^{n+1} = \begin{cases} y'_i & \text{if } i \in Y' \\ \widehat{x}_i & \text{if } i \in I \\ y_i & \text{if } i \in Y \end{cases} \qquad (28)$$

Go to **Step 6.**

    **Step 5.** *(Computing the new trial point if $R_i \leq 0$)*
If $\phi_i(y'_i) \leq 0$ then go to **Step 5.1.** Otherwise go to **Step 5.2.**
    **Step 5.1.** *(Choosing the new trial point if $\phi_i(y'_i) \leq 0$)*
Calculate

$$x^{n+1} = x_{i-1} + \frac{1}{m_i}(z'_{i-1} + \sqrt{z'^2_{i-1} + 2m_i z_{i-1}}), \qquad (29)$$

i.e. the right root of the equation

$$z_{i-1} + z'_{i-1}(x - x_{i-1}) - 0.5 m_i(x - x_{i-1})^2 = 0$$

obtained from the first line of (11) (see Fig.4) and go to **Step 6.**

    **Step 5.2.** *(Choosing the new trial point if $\phi_i(y'_i) > 0$ and $\phi'_i(y'_i)\phi'_i(y_i) < 0$)*
If $\phi'_i(y'_i) \cdot \phi'_i(y_i) \geq 0$ then go to **Step 5.3.** Otherwise if $\phi_i(\widehat{x}_i) > 0$ then compute

$$x^{n+1} = x_i + \frac{1}{m_i}(z'_i + \sqrt{z'^2_i + 2m_i z_i}) \qquad (30)$$

i.e. the right root of the equation

$$z_i - z'_i(x_i - x) - 0.5 m_i(x_i - x)^2 = 0$$

obtained from the third line of (11) (see Fig.5a) and go to **Step 6.**

If $\phi_i(\widehat{x}_i) \leq 0$ then $x^{n+1}$ is calculated following the formula (this situation is presented in Fig.6a)

$$x^{n+1} = \frac{-b_i - \sqrt{b_i^2 - 2m_i c_i}}{m_i} \qquad (31)$$

obtained from the second line of (11) as the left root of the equation

$$0.5 m_i x^2 + b_i x + c_i = 0$$

then go to **Step 6.**

    **Step 5.3.** *(Choosing the new trial point if $\phi'_i(y'_i) \cdot \phi'_i(y_i) \geq 0$)*
If $\phi_i(y_i) > 0$ then calculate $x^{n+1}$ using (30) (see Fig.5b) and go to **Step 6.** Otherwise use (31) for calculating $x^{n+1}$ (see Fig.6b).

    **Step 6.** *(The stopping rule)*
If the stopping rule $|x_i - x_{i-1}| \leq \sigma$, where $\sigma$ is from (4), is fulfilled then **Stop.** Otherwise calculate the value $f(x^{n+1})$ and go to **Step 7** setting $b^{n+1} = x^{n+1}$ if $f(x^{n+1}) < 0$.

**Step 7.** *(Adjusting the search information)*
Calculate the value $f'(x^{n+1})$. Set $n = n + 1$ and go to **Step 1**.
After fulfilment of the stopping rule the following situations can take place:

i. $b^{n+1} \neq b$. This means that we can take $x_\sigma = x_{k-1}$ because this is the last trial point such that $f(x_{k-1}) > 0$.

ii. $b^{n+1} = b$ and $R_i > 0$, for all $i, 1 < i \leq k$. This means that no root has been found in the interval $[a, b]$ and we can continue our investigation taking a new interval $[a^1, b^1]$, where $a^1 = b$. The point

$$x_\sigma^n = argmin\{f(x_j) : 1 \leq j \leq n\}$$

can be taken as a $\sigma$–aproximation of the global minimizer $x'$ over $[a, b]$ and the value $f(x_\sigma^n)$ can be used as an estimate of reliability of our device over the interval $[a, b]$.

iii. $b^{n+1} = b$ and there exists an interval $j$ such that its characteristic $R_j \leq 0$. This situation means that it is necessary to take new $\sigma^1 < \sigma$ because the algorithm stops within the interval $[x_{j-1}, x_j]$ with properties $z_{j-1} > 0$, $z_j > 0$, $R_j \leq 0$ and cannot proceed because $x_{j-1} - x_j \leq \sigma$. For a better understanding of the algorithms logic we present their flow-chart in Fig.7.

Let us say a few words about parameters of the second method. The parameter $\xi$ is a small number reflecting our supposition that $f'(x)$ is not constant over the interval $[x_{i-1}, x_i]$ and $r > 1$ is a reliability parameter of the method. Increasing $r$ means that we suppose that the information obtained during the search is not sufficiently reliable and the objective function behaviour is worse than is seen from the search. Our experience shows that taking $r \in [1.2, 2]$ gives good results both in terms of convergence and in terms of speed.

## 4 Convergence analysis

In this section we demonstrate that the infinite trial sequence $\{x^n\}$ generated by A1 or A2 in the case $\sigma = 0$ has as its limit points (points of accumulation):

-the point $x^*$ from (2) if within $[a, b]$ there exists at least one root of the equation $f(x) = 0$;

-the global minimizer $x'$ if (5) takes place.

We start with establishing these results for A1.

**Theorem 1.** *If there exists the point $x^* \in [a, b]$ from (2), then $x^*$ will be the unique limit point of the sequence $\{x^n\}$ of trial points generated by A1.*

**Proof** Since due to the A1 scheme for all $n \geq 1$ we use

$$m_{i(n)} = K, L \leq K < \infty, 1 < i \leq n,$$

then the auxiliary functions $\phi_i(x)$ from (11) constructed by the algorithm will be the support ones for all $i, 1 < i \leq k(n)$, where k is from (9), i.e.

$$f(x) \geq \phi_i(x), x \in [x_{i-1}, x_i], 1 < i \leq k(n).$$

Denote by $t = t(n)$ the number of an interval $[x_{t-1}, x_t]$ such that

$$x^* \in [x_{t-1}, x_t] \tag{32}$$

in the course of the $n$th iteration. Due to (16), (17) its characteristic $R_t$ is such that

$$R_t < 0 \tag{33}$$

Since $K < \infty$ there exists an infinite sequence of iteration numbers $\{d\}$ such that

$$R_j > 0, 1 < j < t(d), d \in \{d\}. \tag{34}$$

This means that every trial with the number $d + 1, d \in \{d\}$, will fall in the interval $[x_{t-1}, x_t]$. Using again the inequality $L \leq K < \infty$ and (29)–(31) we obtain that

$$\lim_{d \to \infty} x^{d+1} = x^*.$$

To conclude the proof we show that $x^*$ is the unique accumulation point of $\{x^n\}$. Presence of another limit point $\bar{x}$ on the right of $x^*$ is impossible because of **Step 1** of A1 and (9). The situation $\bar{x} < x^*$ cannot take place for the following reason.

Let $\bar{x} \in [x_{c(n)-1}, x_{c(n)}]$ in the course of the $n$th iteration. Then, if $\bar{x}$ is a limit point of $\{x^n\}$ the characteristic $R_{c(n)}$ of the interval $[x_{c(n)-1}, x_{c(n)}]$ should be less than 0 infinitely many times, but this is impossible because of (2) and limitness of $K$. □

**Theorem 2.** *If (5) takes place then all global minimizers will be limit points of the trial sequence $\{x^k\}$ generated by A1.*

**Proof** Due to (5) and the fact that the constant $K < \infty$ there is an iteration number $p$ such that

$$R_j > 0, 1 < j \leq p. \tag{35}$$

This means that for $n > p$ **Step 5** will never be executed and A1 functions as the global optimization method 1 from [21], where the corresponding convergence results are given. □

Let us consider now the performance of the algorithm A2. First of all note that for a correct functioning of the method it is necessary to choose $m_i$ in accordance with the information obtained from the trials executed at the points $x_{i-1}, x_i$. If $m_i$ is underestimated it is possible to obtain $y_i, y_i' \notin [x_{i-1}, x_i]$.

**Proposition 1.** *The choice of $m_i$ by the formula (22) provides that $y_i, y_i' \in [x_{i-1}, x_i]$.*

**Proof** The accordance of the choice of $m_i$ with the local information is done by presence of $\lambda_i$ in (22). This value is determined by (23), (24). A complete discussion of this result can be found in [21]. $\square$

**Theorem 3.** *Let $L_t$ be the local Lipschitz constant of $f(x)$ over the interval $[x_{t-1}, x_t] \ni x^*$, $t = t(n)$, during the nth iteration of A2. If there exists an iteration number $n^*$ such that for all $n > n^*$ the inequality*

$$m_t \geq L_t \tag{36}$$

*holds then the point $x^*$ will be the unique limit point of the trial sequence $\{x^n\}$ generated by A2.*

**Proof** As the values $m_j, 2 \leq j < t$, are bounded (see (22)) as follows:

$$r\xi \leq m_j \leq r \cdot max\{\xi, L\}, 2 \leq j < t, \tag{37}$$

then there exists an iteration number $\overline{n}$ after which a sequence $\{d\}$ from (34) will exist and (34) will take place. Thus, considering iterations with numbers $n > \{n^*, \overline{n}\}$ we obtain that both (33) and (34) hold and the theorem is proved following the proof of theorem 4.1. $\square$

**Remark 1.** *Note that to have convergence to the point $x^*$ it is not necessary to estimate the global Lipschitz constant correctly over the whole region $[a,b]$. It is enough to do it only for the local constant $L_i$ for the subinterval $[x_{t-1}, x_t]$. This condition is significantly weaker than the corresponding convergence results for the methods using estimates of Lipschitz constants (see [8], [9], [16], [22]).*

**Theorem 4.** *Let (5) take places and $L_t$ is the local Lipschitz constant of $f(x)$ over the interval $[x_{t-1}, x_t] \ni x'$, where $x'$ is a global minimizer and there exists a number $n'$ such that (36) takes place. Then $x'$ will be the limit point of the trial sequence generated by A2.*

**Proof** It follows from (37) and (5) that there exists a number $p$ from (35) such that (35) holds. From (36) we obtain

$$\phi_t(x) \leq f(x), x \in [x_{t-1}, x_t].$$

Thus, from the iteration number $\hat{x} = max\{p, n'\}$ A2 functions as the global optimization algorithm 3 from [21] and (36) is its sufficient condition of convergence to the point $x'$. $\square$

## 5 Numerical experiments

In this section we present the results of numerical experiments carried out in order to demonstrate the performance of the new algorithms and to compare them with the grid technique mainly used by engineers to solve the problem (1)–(4).

In the first series of experiments 20 test functions were chosen over the interval [0.2, 7]. Their analytic expressions and characteristics are given in Tab.I. We denote by FRL the first root from the left. The functions are reported in Tab.II and concern general real cases which can be found in many different applications such as filtering and harmonic analysis in electrical or electronic systems, image processing, Wavelet theory, and so on (see [1], [3], [5], [6], [11]–[14], [18], [23]).

The parameters of the algorithms have been chosen as follows: $\xi = 10^{-6}$, $r = 1.2$ for the algorithm A2 and $\sigma = 10^{-4}(b-a)$ for the algorithms A1 and A2 and for the grid method. We used exact Lipschitz constants for $f'(x)$ in A1 in all the experiments. Tab.II contains the numbers of trials required by A1, A2 and the grid method working with the step $\sigma$ before satisfaction of the stopping rule.

In the second part of the experiments we solved practical electrotechnical problems by finding the cutoff frequency for the filters presented in Section 2. The parameters for the Chebyshev filter were the following: $R = 1\Omega; L = 2H; C = 4F$. The cutoff frequency was found as the first zero crossing for the function:

$$f(\omega) = F(\omega)^2 - 0.5F_{max}^2 \tag{38}$$

(see Fig.8), where $F(\omega)$ is from (7) and it was found in the point $\omega = 0.8459 rad/s$. This result was obtained in 2745 iterations by the grid method, in 11 iterations by the algorithm A1 and in 10 iterations by the algorithm A2.

The second filter is a passband filter (see Section 2). The parameters for this filter have been chosen as follows: $R_1 = 3108\Omega$, $L_1 = 40e^{-3}H$, $C_1 = 1e^{-6}F$, $R_2 = 477\Omega$, $L_2 = 350e^{-2}H$, $C_2 = 0.1e^{-6}F$. The cutoff frequency was found as the first zero crossing for the function:

$$f(\omega) = -(F(\omega)^2 - 0.5F_{max}^2) \tag{39}$$

(see Fig.9), where $F(\omega)$ is from (8) and it was found in the point $\omega = 4824.43 rad/s$. This result was obtained in 4474 iterations by the grid method, in 44 iterations by the algorithm A1 and in 27 iterations by the algorithm A2.

## 6  Conclusions

In this paper we have considered the problem of finding the first root from the left of an equation $f(x) = 0$, where $f(x)$ satisfies condition (1). This problem very often arises in practice and we have presented two applications from signal filtering.

For solving this problem we have proposed two methods. The first one uses the exact *a priori* given Lipschitz constant $L$. When $L$ is not known a priori the second method solves the problem using adaptive estimation of the local Lipschitz constant in the course of the search. It uses the obtained estimates in order to accelerate the search.

Numerical experiments executed with real problems and with a set of test functions demonstrate good performance of the new techniques in comparison with the method usually used by engineers. Comparing numerically the first algorithm with the second we can see that, on the set of functions considered in the experiments, the use of local estimates accelerates the search.

**Acknowledgment.** This research was partially supported by the Russian Foundation of Fundamental Research grant 95-01-01073. We would like to thank the first referee for the helpful suggestions.

# References

[1] G. ANTONELLI, F. BINASCO, G. DANESE, D. DOTTI, *Virtually Zero Cross-Talk Dual Frequency Eddy Current Analyzer Based on Personal Computer*, IEEE Transaction on Instrumentation and Measurement, 43 (1994), pp. 463–468.

[2] W. BARITOMPA, *Accelerations for a Variety of Global Optimization Methods*, J. of Global Optimization, 4(1) (1994), pp. 37–45.

[3] D. BEDROSIAN, J. VLACH, *Time Domain Analysis of Network with Internally Controlled Switches*, IEEE Trans. Power Electron, PE–5(3) (1990).

[4] L. BREIMAN, A. CUTLER, *A Deterministic Algorithm for Global Optimization*, Math. Programming, 58 (1993), pp. 179–199.

[5] L. O. CHUA, Charles A. DESOER, Ernest S. KUH, *Linear and Non linear Circuits*, MacGraw Hill, Singapore, 1987.

[6] L. D. COSART, L. PEREGRINO, A. TAMBE, *Time Domain Analysis and its Practical Application to the Measurement of Phase Noise and Jitter*, IEEE Instrumentation and Measurement Technology Conference, Brussels, Belgium, (1996), pp. 430–1435.

[7] C.A. FLOUDAS, P.M. PARDALOS, *State of the Art in Global Optimization*, Kluwer Academic Publishers, Dordrecht, 1996.

[8] E.A. GALPERIN, *The Alpha Algorithm and the Application of the Cubic Algorithm in Case of Unknown Lipschitz Constant*, Computers Math. Applic., 25(11–12) (1993), pp. 71–78.

[9] V.P. GERGEL, *A Global Search Algorithm Using Derivatives*, Systems Dynamics and Optimization, N. Novgorod University Press, 1992, pp.161–178.

[10] R. HORST, P.M. PARDALOS, *Handbook of Global Optimization*, Kluwer Academic Publishers, Dordrecht, 1995.

[11] D. E. JOHNSON, *Introduction to Filter Theory*, Prentice Hall Inc., New Jersey, 1976.

[12] H.Y-F. LAM, *Analog and Digital Filters-Design and Realization*, Prentice Hall Inc., New Jersey, 1979.

[13] A.M. LUCIANO, A.G.M. STROLLO, *A Fast Time-Domain Algorithm for the Simulation of Switching Power Converters*, IEEE Trans. on Power Electr., PE–5 (3) (1990), pp. 363–370.

[14] S. MALLAT, *Zero-Crossing of a Wavelet Transform*, IEEE Trans. on Inf. Theory, 37(4) (1991), pp. 1019–1033

[15] S.A. PIJAVSKII, *An Algorithm for Finding the Absolute Extremum of a Function*, USSR Math. Math. Physics, 12, 1972, pp.57–67.

[16] J. PINTER, *Global Optimization in Action*, Kluwer Academic Publisher, 1996.

[17] W.H. PRESS, Y.B.P. FLANNER, S.A. TEUKOLSKY, W.T. VETTERLING, *Numerical Recipes: the Art of Scientific Computing*, Cambridge University Press, Cambridge, 1986.

[18] —, *Proc. of IEE Colloquium on Transputer Application*, IEE, London, UK, 1989.

[19] Ya.D. SERGEYEV, *A One-Dimensional Deterministic Global Minimization Algorithm*, Comput. Maths. Math. Phys, 35(5) (1995), pp. 705–717.

[20] Ya.D. SERGEYEV, *An Information Global Optimization Algorithm with Local Tuning*, SIAM J. Opt., 5(4) (1995), pp. 858–870.

[21] Ya.D. SERGEYEV, *Global One-Dimensional Optimization Using Smooth Auxiliary Functions*, to appear in Mathematical Programming, (1998).

[22] R.G. STRONGIN, *Numerical Methods on Multiextremal Problems*, Nauka, Moscow, 1978.

[23] P. TURCZA, R. SROKA, T. ZIELINSKI, *Implementation of an Analytic Signal Method of Instantaneous Phase Detection in Real-Time on Digital Signal Processor*, TC-4 IMEKO Modern Electrical and Magnetic Measurement, Prague, 8 (1995), pp. 482–486.

[24] G.R. WOOD, B.P. ZHANG, *Estimation of the Lipschitz Constant of a Function*, J. of Global Optimization, 8 (1996), pp. 91–103.

[25] Ya.D. SERGEYEV, *A Global Optimization Algorithm Using Derivatives and Local Tuning*, ISI–CNR Report, 1 (1994), Rende, Italy.

[26] Ya.D. SERGEYEV, *Global Optimization Algorithms Using Smooth Auxiliary Functions*, ISI–CNR Report, 5 (1994), Rende, Italy.

[27] Ya.D. SERGEYEV, *A Method Using Local Tuning for Minimizing Functions with Lipschitz Derivatives*, Developments in Global Optimization, eds. E. Bomze, T. Csendes, R. Horst and P.M. Pardalos, Kluwer Academic Publishers, (1994), pp. 199–216.

[28] D. MacLAGAN, T. STURGE, W. BARITOMPA, *Equivalent Methods for Global Optimization*, State of Art in Global Optimization, eds. C.A. Floudas, P.M. Pardalos, (1996), pp. 201–212.

List of figures

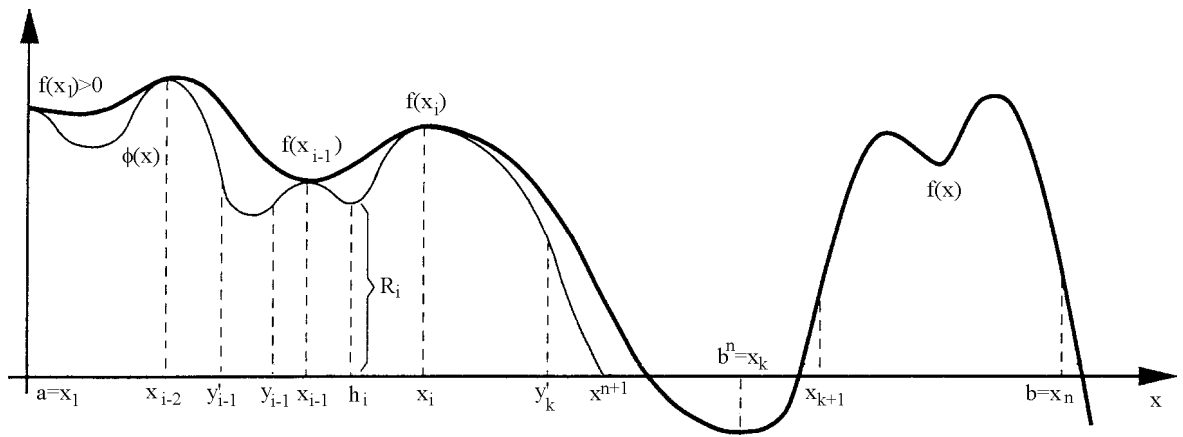

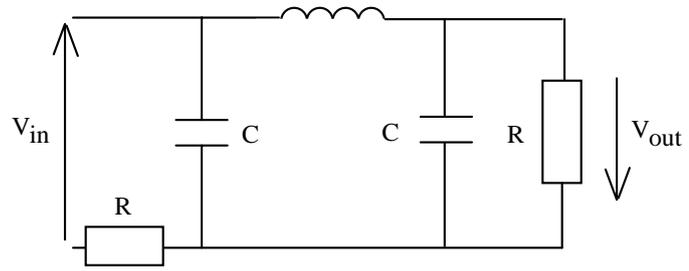

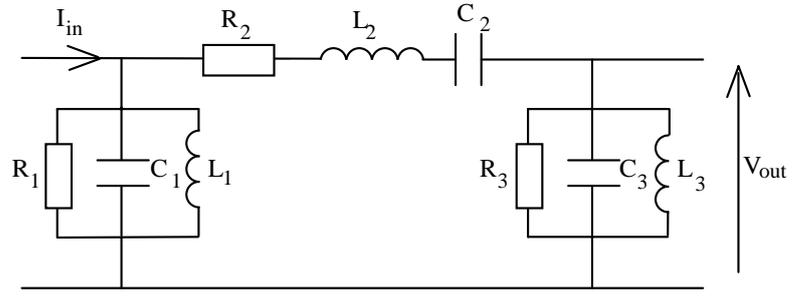
Fig.2

Fig.3

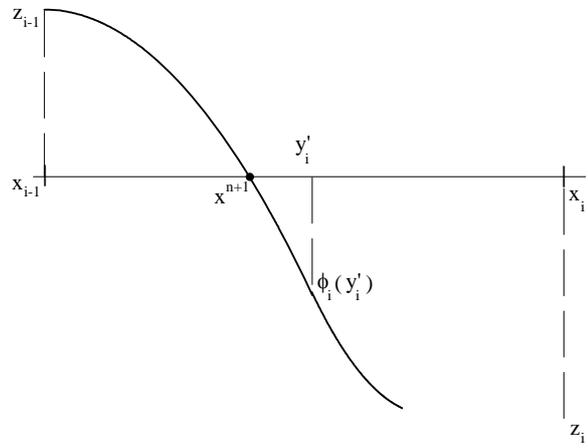

Fig.4

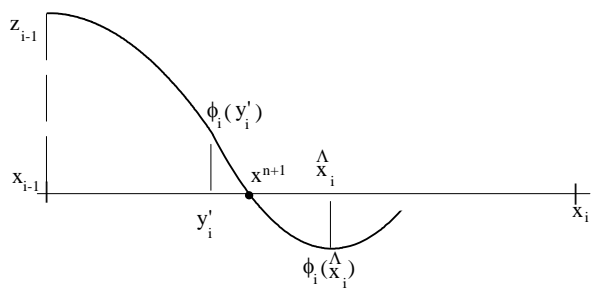

a)

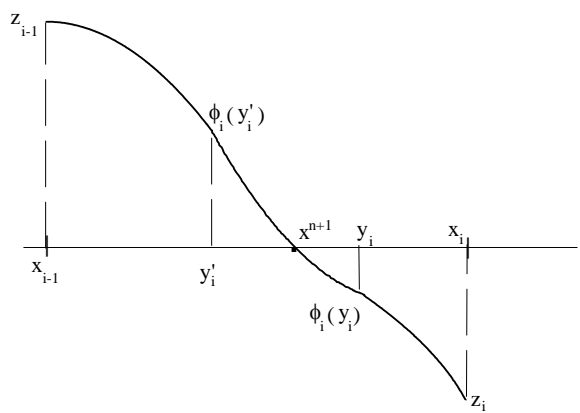

b)

Fig.5

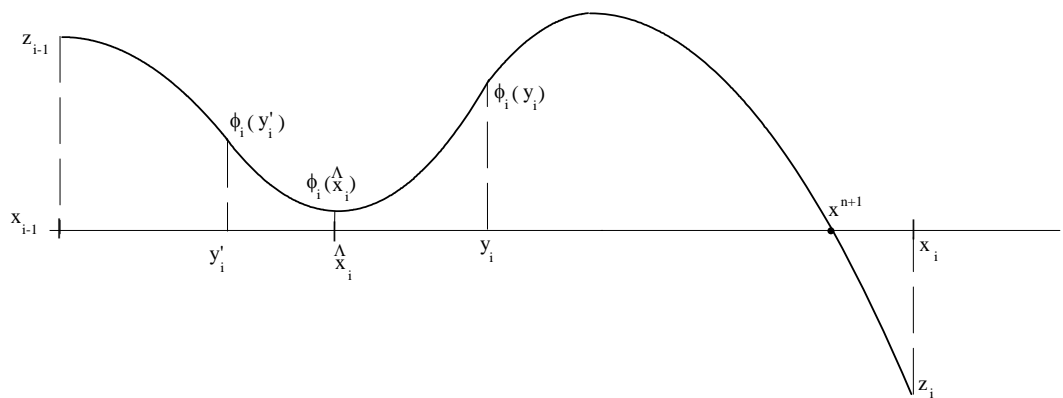

a)

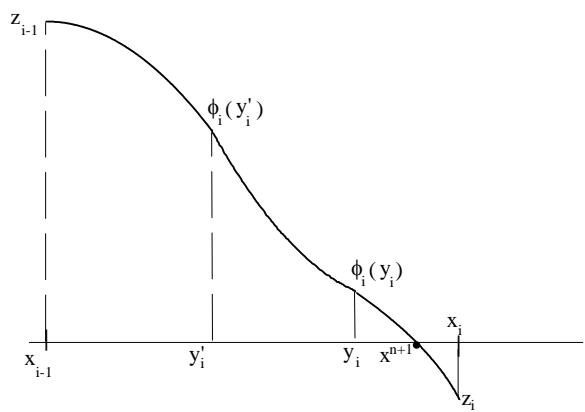

b)

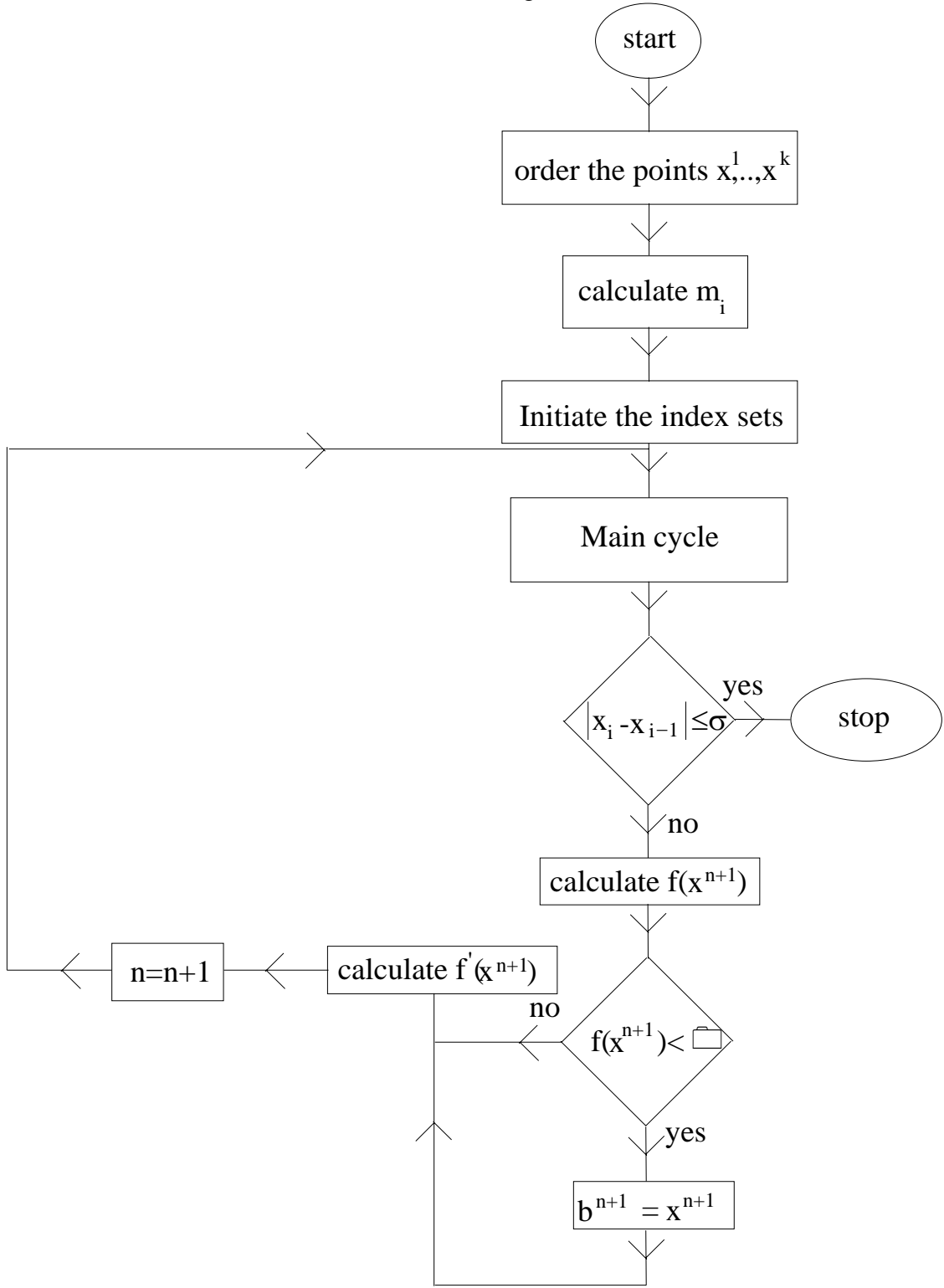

Fig.6

Fig.7 a)

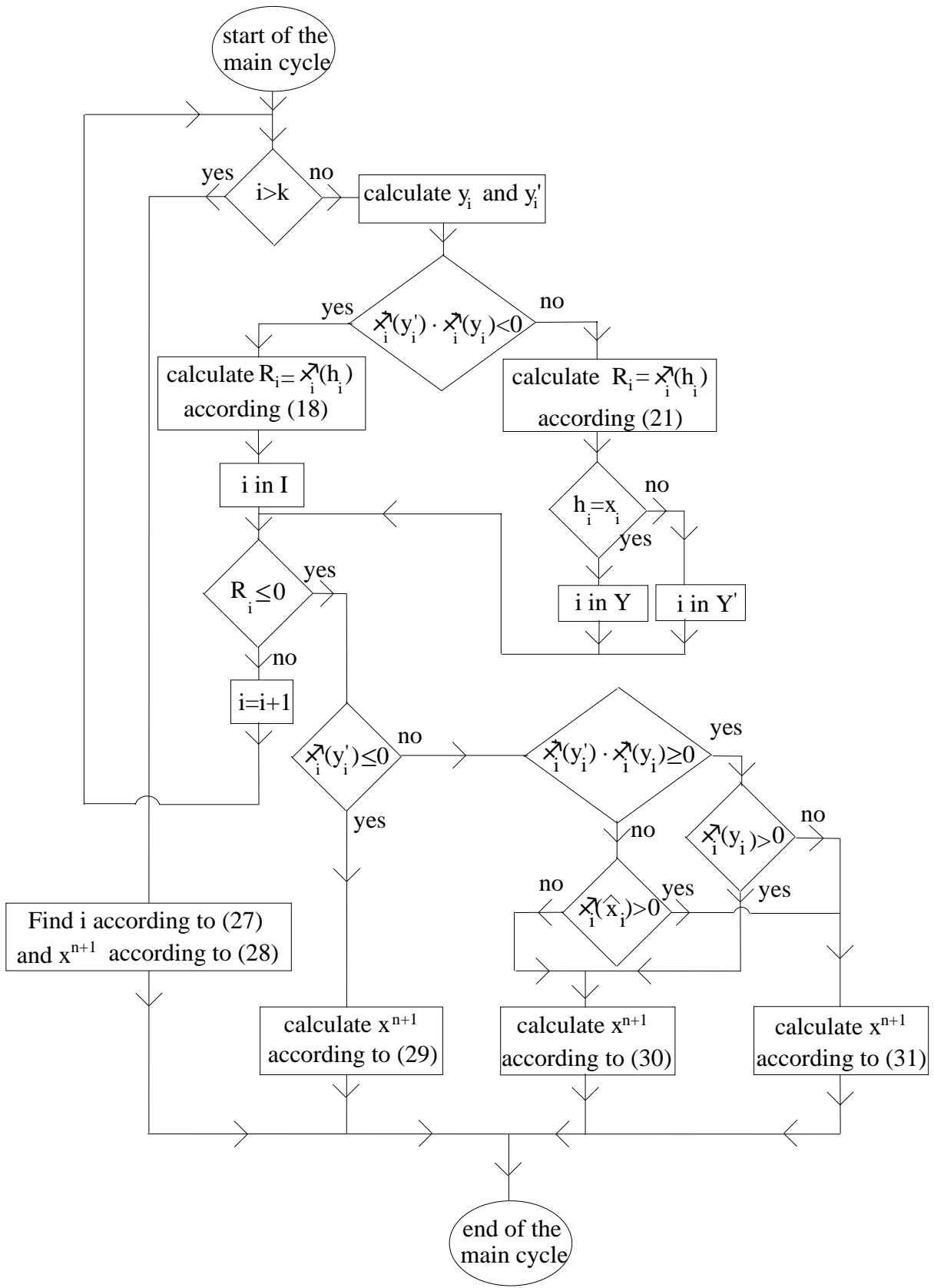

Fig. 7 b)

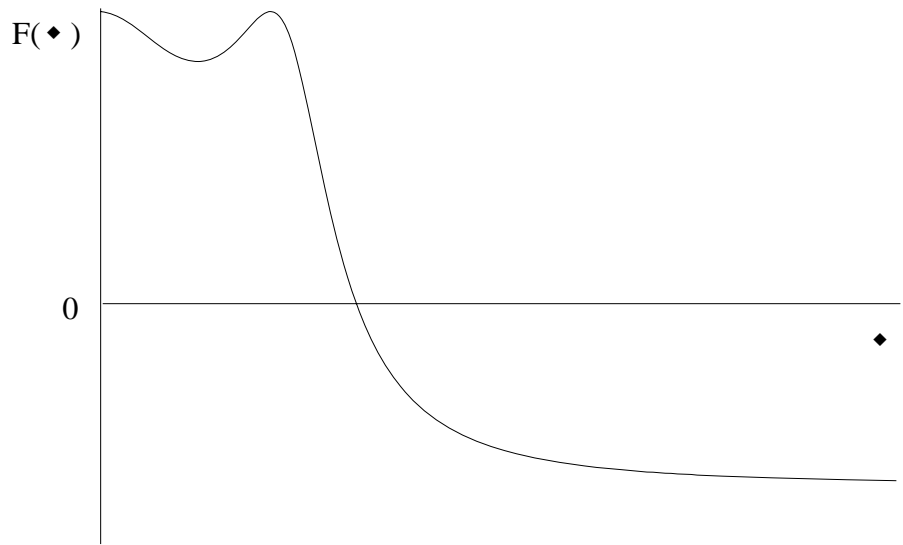

Fig.8

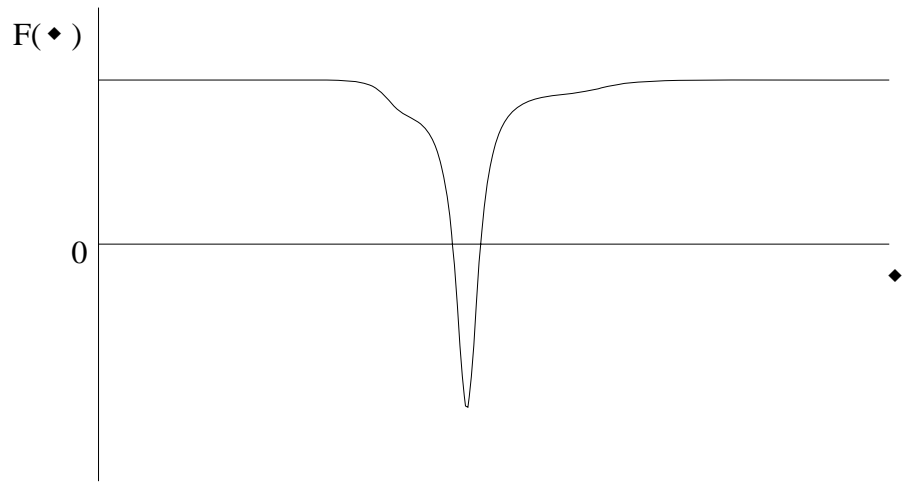

Fig.9

List of Tables

| N° | Function f(x) | Number of roots | FRL | Number of local extrema |
|---|---|---|---|---|
| 1 | $-0.5x^2 \ln(x) + 5$ | 1 | 3.0117 | 3 |
| 2 | $-e^{-x}\sin(2\pi x) + 1$ | - | - | 13 |
| 3 | $-\sqrt{x} \cdot \sin(x) + 1$ | 3 | 1.17479 | 4 |
| 4 | $x\sin(x) + \sin(10x/3) + \ln(x) - 0.84x + 1.3$ | 2 | 2.96091 | 6 |
| 5 | $x + \sin(5x)$ | 2 | 0.82092 | 13 |
| 6 | $-x \cdot \sin(x) + 5$ | - | - | 4 |
| 7 | $\sin(x)\cos(x) - 1.5\sin^2(x) + 1.2$ | 4 | 1.34075 | 7 |
| 8 | $2\cos(x) + \cos(2x) + 5$ | - | - | 6 |
| 9 | $2 \cdot \sin(x) \cdot e^{-x}$ | 2 | 3.1416 | 4 |
| 10 | $(3x - 1.4)\sin(18x) + 1.7$ | 34 | 1.26554 | 42 |
| 11 | $(x+1)^3 / x^2 - 7.1$ | 2 | 1.36465 | 3 |
| 12 | $\begin{cases} \sin(5x) + 2 & x \leq \pi \\ 5\sin(x) + 2 & x > \pi \end{cases}$ | 2 | 3.55311 | 8 |
| 13 | $e^{\sin(3x)}$ | - | - | 9 |
| 14 | $\sum_{k=0}^{5} k\cos[(k+1)x + k] + 12$ | 2 | 4.78308 | 15 |
| 15 | $2(x-3)^2 - e^{x/2} + 5$ | 2 | 3.281119 | 4 |
| 16 | $-e^{\sin(x)} + 4$ | - | - | 4 |
| 17 | $\sqrt{x}\sin^2(x)$ | 4 | 3.141128 | 6 |
| 18 | $\cos(x) - \sin(5x) + 1$ | 6 | 1.57079 | 13 |
| 19 | $-x - \sin(3x) + 1.6$ | 3 | 1.96857 | 9 |
| 20 | $\cos(x) + 2\cos(2x)e^{-x}$ | 2 | 1.14071 | 4 |

Tab.I

| N° | Function | Grid | A1 | A2 |
|---|---|---|---|---|
| 1 | 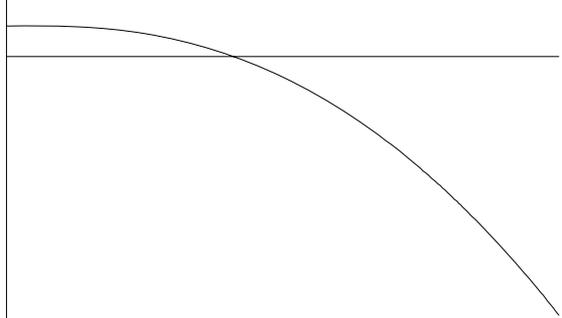 | 4135 | 5 | 5 |
| 2 | 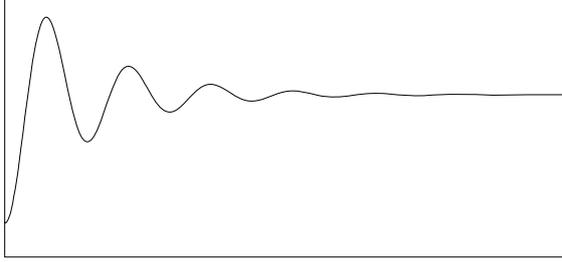 | 10000 | 31 | 34 |
| 3 | 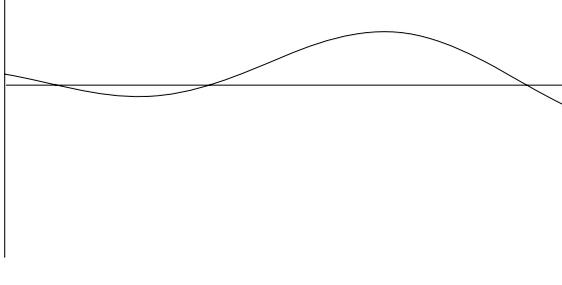 | 1295 | 6 | 5 |
| 4 | 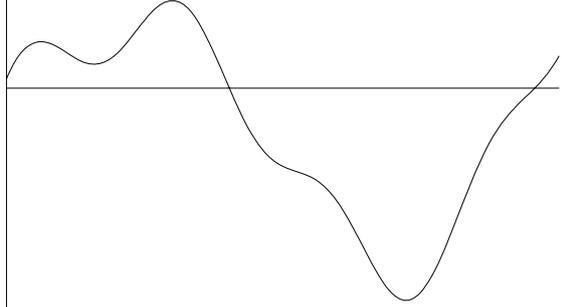 | 4060 | 12 | 7 |
| 5 | 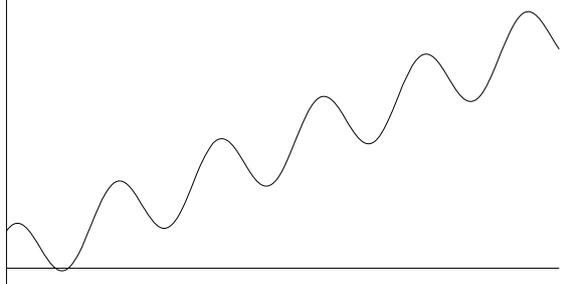 | 5470 | 7 | 11 |

Tab.IIa Two methods for solving optimization problems arising in electronic measurements and electrical engineering

| N | Function | Grid | A1 | A2 |
|---|---|---|---|---|
| 6 | 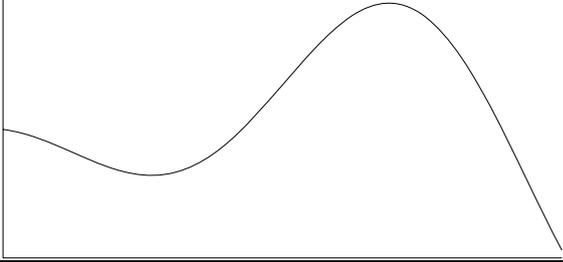 | 10000 | 10 | 9 |
| 7 | 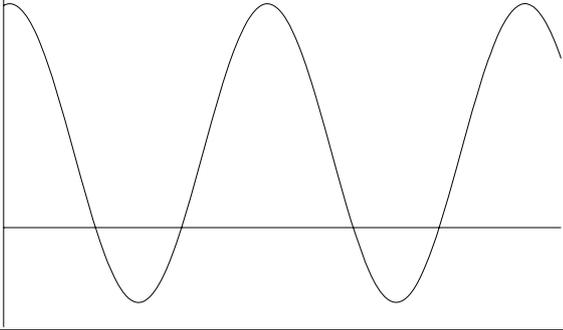 | 1678 | 5 | 6 |
| 8 | 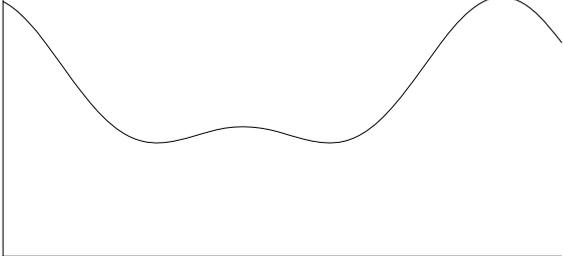 | 10000 | 36 | 24 |
| 9 | 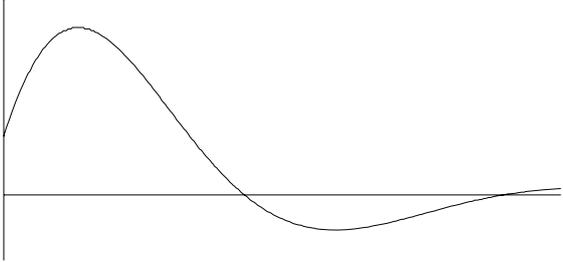 | 4326 | 15 | 10 |
| 10 | 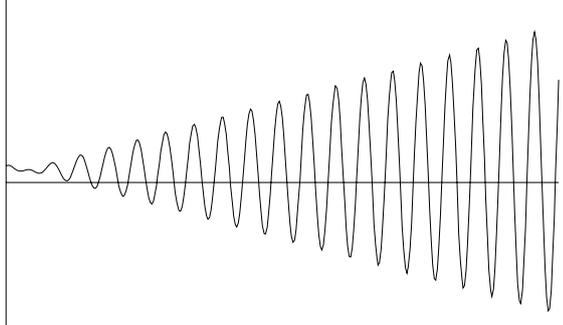 | 1567 | 55 | 12 |

Tab.IIb Two methods for solving optimization problems arising in electronic measurements and electrical engineering

| N | Function | Grid | A1 | A2 |
|---|---|---|---|---|
| 11 | 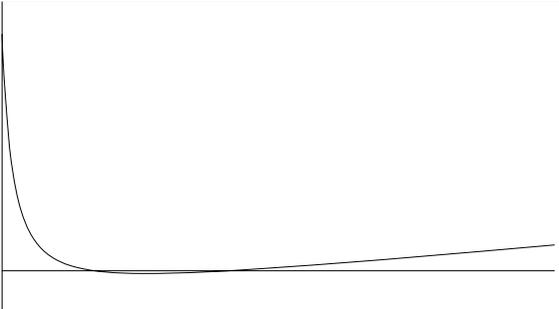 | 1713 | 69 | 60 |
| 12 | 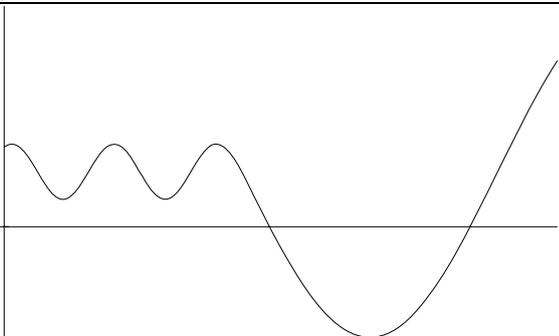 | 4931 | 13 | 6 |
| 13 | 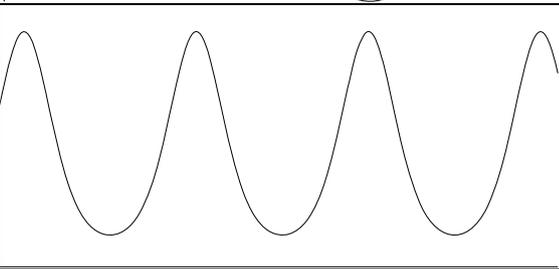 | 10000 | 99 | 39 |
| 14 | 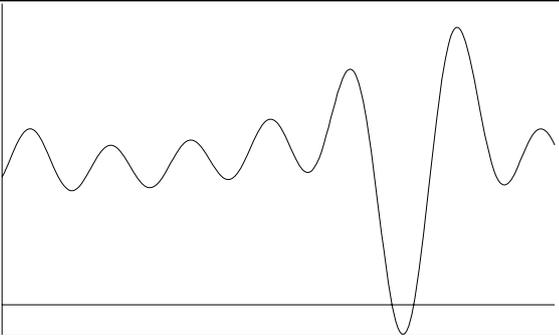 | 6740 | 23 | 18 |
| 15 | 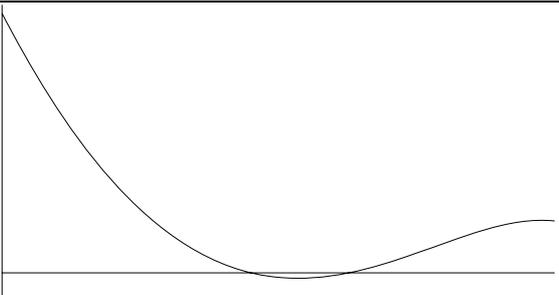 | 4531 | 9 | 9 |

Tab.IIc Two methods for solving optimization problems arising in electronic measurements and electrical engineering

| N | Function | Grid | A1 | A2 |
|---|---|---|---|---|
| 16 | 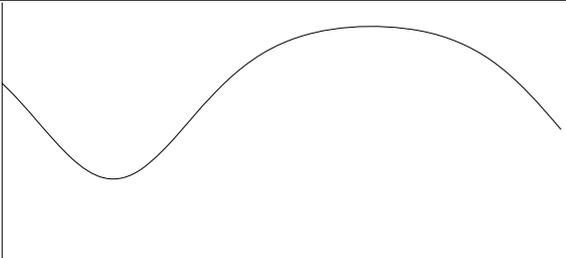 | 10000 | 7 | 12 |
| 17 | 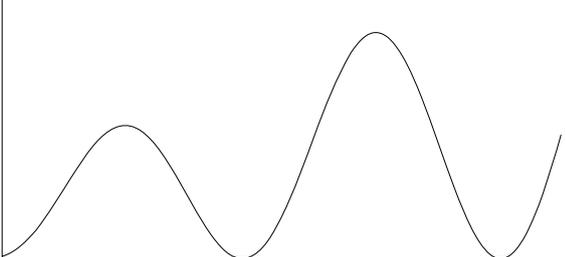 | 4325 | 20 | 17 |
| 18 | 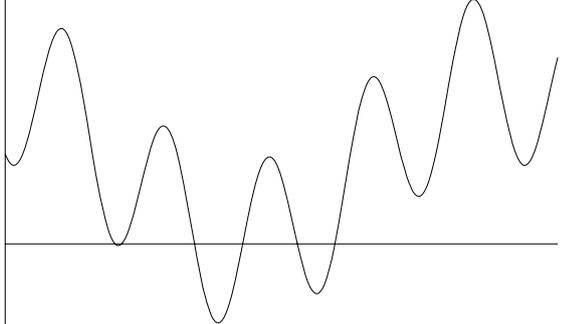 | 2016 | 11 | 10 |
| 19 | 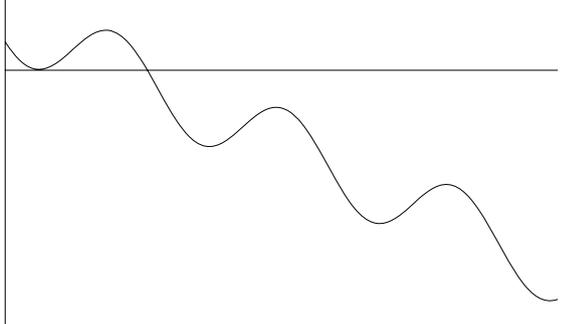 | 2601 | 12 | 12 |
| 20 | 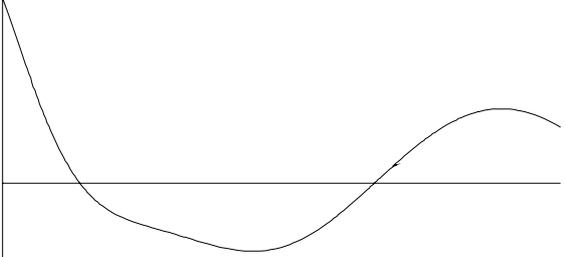 | 7413 | 6 | 6 |
| Average | | 5119.15 | 22.55 | 16.17 |

Tab.IId Two methods for solving optimization problems arising in electronic measurements and electrical engineering